\documentclass[11pt]{article}
\usepackage [latin1]{inputenc}
\usepackage{amssymb}
\usepackage{amsfonts}
\usepackage{amsmath}

\newcommand{\C}{\mathbb{C}}

\newcommand{\Z}{\mathbb{Z}}

\newcommand{\Gl}{\mathop{\mathrm{Gl}}\nolimits}
\newcommand{\gl}{\mathop{\mathrm{gl}}\nolimits}
\newcommand{\spp}{\mathop{\mathrm{sp}}\nolimits}

\newcommand{\Spp}{\mathop{\mathrm{Sp}}\nolimits}

\newcommand{\Ad }{\mathop{\mathrm{Ad}}\nolimits}

\newcommand{\dee}{\mathop{\! \, \mathrm{d}\!}\nolimits}

\newcommand{\comp}{\, \raisebox{2pt}{$\scriptstyle\circ \, $}}
\newcommand{\setrule}{\, \rule[-4pt]{.5pt}{13pt}\, }
\newcommand{\onehalf}{\mbox{$\frac{\scriptstyle 1}{\scriptstyle 2}\,$}}

\newcommand{\lefthook}{\mbox{$\, \rule{8pt}{.5pt}\rule{.5pt}{6pt}\, \, $}}

\newcommand{\smalldbydt}{\mbox{${\scriptstyle \frac{\dee}{\dee t}}
\rule[-6pt]{.5pt}{12pt} \raisebox{-6pt}{$ \, {\scriptscriptstyle t=0}$}$}}
\newcommand{\vvee}{\mbox{\tiny $\vee $}}

\begin{document}
\begin{center}
\noindent{\Large \bf A momentum map for the Heisenberg group} \newline 
\rule{0pt}{20pt}{Richard Cushman}\footnotemark 
\end{center}
\footnotetext{Department of Mathematics and Statistics.
University of Calgary, Calgary, Alberta, Canada}
\date{\today} \bigskip

\begin{abstract}
We look at a momentum mapping associated to the Heisenberg group. Following Cushman and van der Kallen 
\cite{cushman-vanderKallen} we algebraically classify the coadjoint 
orbits of the Heisenberg group. We obtain one orbit with a continuous parameter called a modulus. We show that the cocycle associated its momentum mapping is the value of a modulus of the coadjoint orbit. We give a representation theoretic description of this modulus.
\end{abstract}

\section{The affine group ${\mathbb{R} }^{2n}$}

Let $({\mathbb{R} }^{2n}, \omega )$ be a real symplectic vector space, where $\omega $ is the 
standard symplectic form on ${\mathbb{R} }^{2n}$ with matrix $J_{2n}=${\tiny $\begin{pmatrix} 
0 & -I_n \\ I_n & 0 \end{pmatrix}$} with respect to the standard symplectic basis 
${\mathfrak{e}}_n = \{ e_1, \ldots , e_n , f_1, \ldots , f_n \}$. Consider the affine action 
\begin{equation}
\Phi : {\mathbb{R} }^{2n} \times ({\mathbb{R} }^{2n}, \omega ) \rightarrow ({\mathbb{R} }^{2n}, \omega ): (a,v) \mapsto v +a. 
\label{eq-one}
\end{equation}
For each $x$ in the abelian Lie algebra $({\mathbb{R} }^{2n}, [ \, \, , \, \, ])$ of the abelian Lie group 
$({\mathbb{R} }^{2n}, +)$ the infinitesimal generator of the 
action $\Phi $ in the direction $x$ is the vector field $X^x(v) = x$. \medskip 

For each $x \in ({\mathbb{R} }^{2n}, [ \, \, , \, \, ])$ let $J^x: {\mathbb{R} }^{2n} \rightarrow \mathbb{R} : v \mapsto \omega (x,v)$. Then 
$\dee J^x(v)y = \omega (x,y) = \omega (X^x(v),y)$, that is, $X^x = X_{J^x}$. 
Hence the action $\Phi $ is Hamiltonian with momentum mapping $J: {\mathbb{R} }^{2n} \rightarrow ({\mathbb{R} }^{2n})^{\ast }$, where 
$J(v)x = J^{x}(v)$. For $y, z \in {\mathbb{R} }^{2n}$ 
\begin{displaymath}
\{ J^y, J^z \} (v)= \dee J^y(v)X^z(v) = \omega (X^y(v), X^z(v)) = \omega (y,z); 
\end{displaymath}
while $J^{[y,z]}(v) = \omega ([y,z],v) = 0$, since $[y,z] =0$. So 
\begin{displaymath}
\Sigma (y,z) = \{ J^y J^z \} (v) - J^{[y,z]}(v) = \omega (y,z)
\end{displaymath}
is the $({\mathbb{R} }^{2n},[\, \, , \, \, ])$ cocycle associated to the momentum map $J$. We cannot use Kostant's 
theorem \cite[Thm. 4.5, p.175-6]{kostant} to determine a representation of $({\mathbb{R} }^{2n}, \omega )$ because 
$J$ is not coadjoint equvariant. \medskip 

To continue we find a central extension $\mathfrak{g} \subseteq {\mathbb{R} }^{2n} \times \mathbb{R} $ of the abelian Lie algebra 
$({\mathbb{R} }^{2n}, [\, \, , \, \, ])$ by the cocycle $\Sigma $, which has Lie bracket $[(x,t), (y,s)] = (0, \Sigma (x,y))$. 
A Lie group $G$ with Lie algebra $\mathfrak{g}$ has multiplication 
\begin{displaymath}
(v,r) \cdot (w,s) = \big( v+w, r+s+ \onehalf \omega (v,w) \big) . 
\end{displaymath}

We now construct a model of $G$ as a subgroup of $\Gl ({\mathbb{R} }^{2n+2}, \mathbb{R} )$. Consider the map 
\begin{displaymath}
\rho : G \rightarrow \Gl ({\mathbb{R} }^{2n+2}, \mathbb{R} ): (v,r) \rightarrow 
\mbox{\footnotesize $\begin{pmatrix} 1 & 0 & 0 \\ v & I_{2n} & 0 \\ r & \onehalf {\omega }^{\sharp}(v) & 1 
\end{pmatrix}$,}
\end{displaymath}
where ${\omega }^{\sharp}(v) \in ({\mathbb{R} }^{2n})^{\ast }$ is given by ${\omega }^{\sharp}(v)w = 
\omega (v, w)$ for every $v$, $w \in {\mathbb{R} }^{2n}$. Since 
\begin{align*}
\rho (v+w,r+s + \onehalf \omega (v,w)) & = 
\mbox{\footnotesize $\begin{pmatrix} 1 & 0 & 0 \\ v+w & I_{2n} & 0 \\ 
r+s + \onehalf \omega (v,w) & \onehalf {\omega }^{\sharp}(v+w) & 1 \end{pmatrix} $}\\
&\hspace{-.75in} = \mbox{\footnotesize $\begin{pmatrix} 1 & 0 & 0 \\ v & I_{2n} & 0 \\ r & 
\onehalf {\omega }^{\sharp}(v) & 1 
\end{pmatrix} \, \begin{pmatrix} 1 & 0 & 0 \\ w & I_{2n} & 0 \\ s & \onehalf {\omega }^{\sharp}(w) & 1 
\end{pmatrix}$} = \rho (v,r) \rho (w,s), 
\end{align*}
$\rho $ is a group homomorphism. But $\rho (v,r) = I_{2n+2}$ implies $(v,r) = (0,0)$, which is 
the identity element of $G$. Hence $\rho $ is an isomorphism onto its image. Since $G$ is 
the $2n+1$ dimensional Heisenberg group ${\mathrm{H}}_{2n+1}$, so is $H = \rho (G)$. \medskip 

The Lie algebra $\mathfrak{h}$ of $H$ is $\{ X^{\vvee}=${\tiny $\begin{pmatrix} 
0 & 0 & 0 \\ x & 0 & 0 \\ \xi & \onehalf {\omega }^{\sharp}(x) & 0 \end{pmatrix}$} \hspace{-5pt} \setrule  
\mbox{\hspace{-2pt}$\xi \in \mathbb{R}$ and $x \in {\mathbb{R} }^{2n}$}$\}$ with Lie bracket $[X^{\vvee}, Y^{\vvee}] = X^{\vvee}Y^{\vvee} - 
Y^{\vvee}X^{\vvee} $. The map 
\begin{displaymath}
\mathfrak{g} \rightarrow \mathfrak{h}: (x, \xi ) \mapsto \mbox{\footnotesize $\begin{pmatrix}
0 & 0 & 0 \\ x & 0 & 0 \\ \xi & \onehalf {\omega }^{\sharp}(x) & 0 \end{pmatrix}$}
\end{displaymath}
is a Lie algebra isomorphism, since $[X^{\vvee}, Y^{\vvee}] =${\tiny $\begin{pmatrix}
0 & 0 & 0 \\ 0 & 0 & 0 \\ \omega (v,w) & 0 & 0 \end{pmatrix}$}. \medskip

The action 
\begin{equation}
{\Phi }^{\vvee}: H \times ({\mathbb{R} }^{2n}, \omega ) \rightarrow ({\mathbb{R} }^{2n}, \omega ): 
\left( \mbox{\footnotesize $\begin{pmatrix} 1 & 0 & 0 \\ v & 0 & 0 \\ r & \onehalf {\omega }^{\sharp}(v) & 1 
\end{pmatrix}$}, w \right) \mapsto w + v 
\label{eq-two}
\end{equation}
is Hamiltonian. For each $X^{\vvee}=${\tiny $\begin{pmatrix} 0 & 0 & 0 \\ x & 0 & 0 \\ \xi & 
\onehalf {\omega }^{\sharp}(x) & 
0 \end{pmatrix}$}$\in \mathfrak{h}$ the infinitesimal generator of the action ${\Phi }^{\vvee}$ in the direction 
$X^{\vvee}$ is the vector field $X^{\vvee}(v) = x$. Let $J^{X^{\vvee}}: {\mathbb{R} }^{2n} \rightarrow \mathbb{R} : v \mapsto 
\omega (x, v) + \xi $. Since 
\begin{displaymath}
\dee J^{X^{\vvee}}(v)w = \omega (x,w) =  \omega (X^{\vvee}(v), w), 
\end{displaymath}
we obtain $X^{X^{\vvee}} = X_{J^{X^{\vvee}}}$. So the action ${\Phi }^{\vvee}$ is Hamiltonian. The 
corresponding momentum mapping is $J:{\mathbb{R} }^{2n} \rightarrow {\mathfrak{h}}^{\ast }$, where 
$J(v)X^{\vvee} = J^{X^{\vvee}}(v)$. The mapping $J$ is coadjoint equivariant, that is, 
$J({\Phi }^{\vvee}_g(v)) = {\Ad }^T_{g^{-1}}J(v)$ for every $g\in G$ and every $v\in {\mathbb{R} }^{2n+2}$. We verify this. Since the Heisenberg group $H$ is 
connected we need only show that $\{ J^{X^{\vvee}}, J^{Y^{\vvee}} \} = J^{[X^{\vvee}, Y^{\vvee}]}$, 
which is the infinitesimalization of the coadjoint equivariance condition. We have 
\begin{align*}
\{ J^{X^{\vvee}}, J^{Y^{\vvee}} \} (v) & = L_{X^{Y^{\vvee}}}J^{X^{\vvee}}(v) = \dee J^{X^{\vvee}}(v)X^{Y^{\vvee}}(v) \\
& = \omega (X^{X^{\vvee}}(v), X^{Y^{\vvee}}(v)) = \omega (x,y) = J^{[X^{\vvee}, Y^{\vvee}]}(v). 
\end{align*}

\section{Coadjoint orbits of ${\mathrm{H}}_{2n+1}$}

In this section we classify the coadjoint orbits of the Heisenberg group ${\mathrm{H}}_{2n+1}$. \medskip 

\vspace{-.15in}First we find a subgroup of $\Spp (V, \mathcal{J})$ which has an isotropy group equal to 
${\mathrm{H}}_{2n+1}$. Here $V = \mathbb{R} \times {\mathbb{R} }^{2n} \times \mathbb{R} $. With respect to the basis 
$\mathfrak{e} = \{ e_0; {\mathfrak{e}}_n; f_{n+1} \}$ of $V$ the matrix of the symplectic form 
$\mathcal{J}$ is {\tiny $\begin{pmatrix} 0 & 0 & 1 \\ 0 & J & 0 \\ -1 & 0 & 0 \end{pmatrix}$}, where 
$J$ is the matrix of a symplectic form on ${\mathbb{R} }^{2n}$ with respect to the basis ${\mathfrak{e}}_n$. \medskip 

Let 
\begin{displaymath}
\widetilde{G} = \left\{ \mbox{\footnotesize $\begin{pmatrix} a & 0 & 0 \\ d & I_{2n} & 0 \\ f & g^T & h \end{pmatrix}$} 
\in \Gl ({\mathbb{R} }^{2n+2}, \mathbb{R} ) \setrule \, a,h \in {\mathbb{R} }^{\times }; \, f\in \mathbb{R}; \, d,g \in {\mathbb{R} }^{2n} \right\} . 
\end{displaymath}
Then $\widetilde{G}$ is a group of real linear mappings of ${\mathbb{R} }^{2n+2}$ into itself which sends 
$\{ 0 \} \times {\mathbb{R}} \subseteq {\mathbb{R} }^{2n+1} \times \mathbb{R} $ into itself and is the identity map on 
$\{ 0 \} \times {\mathbb{R} }^{2n} \times \{ 0 \}$. Let $\widehat{G} = \widetilde{G} \cap 
\Spp ({\mathbb{R} }^{2n+2}, \mathcal{J} )$. Then 
\begin{displaymath}
\widehat{G} = \left\{ \mbox{\footnotesize $\begin{pmatrix} a & 0 & 0 \\ d & I_{2n} & 0 \\ f & a^{-1}(Jd)^T & a^{-1} 
\end{pmatrix} $} \in \Spp ({\mathbb{R} }^{2n+2}, \mathcal{J}) \setrule \, a \in {\mathbb{R} }^{\times }; f \in \mathbb{R}; \, d \in {\mathbb{R} }^{2n} \right\} 
\end{displaymath}
is a Lie group with Lie algebra 
\begin{displaymath}
\widehat{\mathfrak{g}} = \left\{ \mbox{\footnotesize $\begin{pmatrix} \eta  & 0 & 0 \\ \widetilde{x} & 0 & 0 \\
\xi & (J\widetilde{x})^T & -\eta  \end{pmatrix}$} \in \spp ({\mathbb{R} }^{2n+2}, \mathcal{J}) \setrule \, 
\eta, \xi \in \mathbb{R}; \ \widetilde{x} \in {\mathbb{R} }^{2n} \right\} .
\end{displaymath}
The isotropy group ${\widehat{G}}_{f_{n+1}}$ of elements of $\widehat{G}$ which leave the vector 
$f_{n+1}$ fixed is 
\begin{displaymath}
{\widehat{G}}_{f_{n+1}} = \left\{ \mbox{\footnotesize $\begin{pmatrix}
1 & 0 & 0 \\ d & I_{2n} & 0 \\ f & (Jd)^T & 1 \end{pmatrix}$} \in \Spp ({\mathbb{R} }^{2n+2}, \mathcal{J}) \setrule \, 
f \in \mathbb{R}; \, d \in {\mathbb{R} }^{2n} \right\} .
\end{displaymath}
Using the basis ${\mathfrak{e}}_n$ of ${\mathbb{R} }^{2n}$, set $J= \onehalf J_{2n}$. We see that ${\widehat{G}}_{f_{n+1}}$ 
is the Heisenberg group ${\mathrm{H}}_{2n+1}$. The Lie algebra of ${\widehat{G}}_{f_{n+1}}$ is 
\begin{displaymath}
{\widehat{\mathfrak{g}}}_{f_{n+1}} = \left\{ \mbox{\footnotesize $\begin{pmatrix}
0 & 0 & 0 \\ \widetilde{x} & 0 & 0 \\ \xi & (J\widetilde{x})^T & 0 \end{pmatrix}$} \in \spp ({\mathbb{R}}^{2n+2}, \mathcal{J}) \setrule 
\xi \in \mathbb{R} ; \, \widetilde{x} \in {\mathbb{R} }^{2n} \right\} . 
\end{displaymath}

We now begin the classification of the coadjoint orbits of ${\widehat{G}}_{f_{n+1}}$. Let 
$(V = \mathbb{R} \times \widetilde{V} \times \mathbb{R} , \mathcal{J})$ be a real symplectic vector space of 
dimension $2n+2$ with a basis $\mathfrak{e} = \{ e_0; {\mathfrak{e}}_n; f_{n+1} \}$ such that 
the matrix of $\mathcal{J}$ with repect to the basis $\mathfrak{e}$ is {\tiny $\begin{pmatrix} 
0 & 0 & 1 \\ 0 & J & 0 \\ -1 & 0 & 0 \end{pmatrix}$}, where $J$ is one half the matrix $J_{2n}$ of the standard 
symplectic form $\widetilde{\omega }$ on $\widetilde{V}$. A \emph{tuple} $(V,Y, f_{n+1}; \mathcal{J})$ 
is a symplectic vector space $(V, \mathcal{J})$ with $Y \in \widehat{\mathfrak{g}}$, the Lie algebra 
of the Lie group $\widehat{G}$, and $f_{n+1}$ is a basis vector in $\mathfrak{e}$. Two tuples 
$(V,Y, f_{n+1}; \mathcal{J})$ and $(V, Y', f_{n+1}; \mathcal{J})$ are \emph{equivalent} if there is 
a bijective real linear mapping $P: V \rightarrow V$ such that 1) $P \in {\widehat{G}}_{f_{n+1}}$ and 
2) there is a vector $w \in V$ such that $Y' = P(Y+L_{w, f_{n+1}})P^{-1}$. Here 
$L_{w,f_{n+1}} = w \otimes f^{\ast}_{n+1} + f_{n+1} \otimes w^{\ast }$ with $w^{\ast}(z) = z^T\mathcal{J}w$ 
for every $z \in V$. \medskip 

\noindent \textbf{Lemma 1} Let $w = w_0e_0 + \widetilde{w} + w_{2n+1}f_{n+1} \in V$. We have 
$L_{w,f_{n+1}} \in \widehat{\mathfrak{g}}$. \medskip 

\noindent \textbf{Proof} We compute the matrix of $L_{w, f_{n+1}}$ with respect to the basis 
$\mathfrak{e}$. 
\begin{align*}
L_{w,f_{n+1}}(e_0) & = (w\otimes f^{\ast }_{n+1})(e_0) + (f_{n+1} \otimes w^{\ast})(e_0) \\
& = (e^T_0 \mathcal{J}f_{n+1})w + (e^T_0 \mathcal{J} w)f_{n+1} \\
& = w + w_{2n+1}f_{n+1} = w_0 e_0 + \widetilde{w} + 2w_{2n+1}f_{n+1}; \\
L_{w, f_{n+1}}({\mathfrak{e}}_n) & = (w \otimes f^{\ast }_{n+1})({\mathfrak{e}}_n) + 
(f_{n+1} \otimes w^{\ast})({\mathfrak{e}}_n) \\
& = (({\mathfrak{e}}_n)^T \mathcal{J}f_{n+1})w + (({\mathfrak{e}}_n)^T\mathcal{J}w) f_{n+1} \\
& = (J\widetilde{w})^Tf_{n+1}; \\
L_{w, f_{n+1}}f_{n+1}& = (w \otimes f^{\ast}_{n+1})(f_{n+1}) + (f_{n+1} \otimes w^{\ast })f_{n+1} \\
& = (f^T_{n+1} \mathcal{J}f_{n+1})w +(f^T_{n+1}\mathcal{J}w)f_{n+1} = - w_0f_{n+1}. 
\end{align*}
Thus the matrix of $L_{w, f_{n+1}}$ with respect to the basis $\mathfrak{e}$ is 
{\tiny $\begin{pmatrix} w_0 & 0 & 0 \\ \widetilde{w} & 0 & 0 \\ 2w_{2n+1} & (J\widetilde{w})^T & - w_0 \end{pmatrix}$}, which lies in $\widehat{\mathfrak{g}}$. \hfill $\square $ \medskip 

\noindent Hence the definition of equivalence makes sense. \medskip

\noindent \textbf{Corollary 1A} $\widehat{\mathfrak{g}} = \{ L_{w, f_{n+1}} \setrule \, w \in V \}$. \medskip 

\noindent \textbf{Proof} From lemma 1 and the definition of $\widehat{\mathfrak{g}}$ it follows that 
$L_{w,f_{n+1}} \in \widehat{\mathfrak{g}}$ for every $w \in V$. Suppose that $X=${\tiny $\begin{pmatrix}
0 & 0 & 0 \\ \widetilde{x} & 0 & 0 \\ \xi & (J\widetilde{x})^T & 0 \end{pmatrix}$}$\in \widehat{\mathfrak{g}}$. Let 
$w = w_0 +\widetilde{x} +\onehalf \xi f_{n+1}$. Then $X=L_{w, f_{n+1}}$. \hfill $\square $ \medskip 

\noindent \textbf{Lemma 2} For every $P \in \Spp (V, \mathcal{J})$ we have $PL_{w,f_{n+1}}P^{-1} = 
L_{Pw, Pf_{n+1}}$. \medskip 

\noindent \textbf{Proof} For every $z \in V$ 
\begin{align}
L_{Pw, Pf_{n+1}}z & = (Pw \otimes (Pf_{n+1})^{\ast })(z) + (Pf_{n+1} \otimes (Pw)^{\ast})(z) \notag \\
&\hspace{-.5in} = ((z^T\mathcal{J})Pf_{n+1})Pw + (z^T\mathcal{J}Pw)Pf_{n+1} \notag \\
&\hspace{-.5in} = ((P^{-1}z)^T\mathcal{J}f_{n+1})Pw + ((P^{-1}z)^T\mathcal{J}w)Pf_{n+1}, \quad 
\mbox{since $P\in \Spp (V, \mathcal{J})$} \notag \\
&\hspace{-.5in} = P(w \otimes f^{\ast }_{n+1} +f_{n+1} \otimes w^{\ast })P^{-1}(z) = (PL_{w, f_{n+1} }P^{-1})z. 
\tag*{$\square $}
\end{align}

Being equivalent is an equivalence relation on the set of tuples. An equivalence class $\nabla $ of a 
set of tuples is a \emph{cotype}. $\nabla $ is \emph{represented by} the tuple $(V,Y, f_{n+1}; \mathcal{J})$ if and 
only if $(V,Y, f_{n+1}; \mathcal{J}) \in \nabla $. If the element $Y$ of the tuple $(V,Y, f_{n+1}; \mathcal{J})$ 
is nilpotent with $Y^{m+1}V=0$ but $Y^m V \ne 0$ for some $m \in {\Z}_{\ge 1}$, then the tuple is 
nilpotent of height $m$. Since equivalent nilpotent tuples have the same height, we say that the 
corresponding cotype is \emph{nilpotent of height $m$}. \medskip 

For each $Y \in \widehat{\mathfrak{g}}$ let ${\ell }_Y: \widehat{\mathfrak{g}} \rightarrow \mathbb{R} : Z \mapsto 
\mathrm{tr}\, YZ$. Then ${\ell }_Y \in {\mathfrak{g}}^{\ast }$. \medskip 

\noindent \textbf{Theorem 3} The correspondence 
\begin{equation}
(V,Y,f_{n+1}; \mathcal{J}) \mapsto ({\ell }_Y)_{|{\widehat{\mathfrak{g}}}_{f_{n+1}} }
\label{eq-three} 
\end{equation}
between tuples and elements of $({\widehat{\mathfrak{g}}}_{f_{n+1}})^{\ast }$ induces a bijection between 
cotypes and ${\widehat{G}}_{f_{n+1}}$ coadjoint orbits. \medskip 

To prove the theorem we need \medskip 

\noindent \textbf{Fact 4} 
\begin{displaymath}
({\widehat{\mathfrak{g}}}_{f_{n+1}})^{\circ} = \{ X \in \widehat{\mathfrak{g}} \setrule \, 
{\ell }_X(Y) =0 \, \mbox{for every $Y=${\tiny $\begin{pmatrix} 0 & 0 & 0 \\ v & 0 & 0 \\ \xi & (Jv)^T & 0 
\end{pmatrix}$}$\in {\widehat{\mathfrak{g}}}_{f_{n+1}} $} \} = \widehat{\mathfrak{g}}. 
\end{displaymath}

\noindent \textbf{Proof} Suppose that $X=${\tiny $\begin{pmatrix} \zeta & 0 & 0 \\ d & 0 & 0 \\ \eta  & (Jd)^T  & -\zeta  
\end{pmatrix}$}$\in \widehat{\mathfrak{g}}$. For every $Y \in {\widehat{\mathfrak{g}}}_{f_{n+1}}$ we have 
\begin{displaymath}
{\ell }_X(Y) = \mathrm{tr}\, XY = \mathrm{tr}\, \left[ \mbox{\footnotesize $\begin{pmatrix} 
\zeta & 0 & 0 \\ d & 0 & 0 \\ \eta & (Jd)^T & -\zeta  \end{pmatrix} \, \begin{pmatrix} 
0 & 0 & 0 \\ v & 0 & 0 \\ \xi & (Jv)^T & 0 \end{pmatrix}$} \right] =0. 
\end{displaymath}
So $X \in ({\widehat{\mathfrak{g}}}_{f_{n+1}})^{\circ}$. Hence $\widehat{\mathfrak{g}} \subseteq 
({\widehat{\mathfrak{g}}}_{f_{n+1}})^{\circ}$. By definition $({\widehat{\mathfrak{g}}}_{f_{n+1}})^{\circ} 
\subseteq \widehat{\mathfrak{g}}$. Thus $({\widehat{\mathfrak{g}}}_{f_{n+1}})^{\circ} = \widehat{\mathfrak{g}}$. 
\hfill $\square $ \medskip

\noindent \textbf{Proof of theorem 3.} Let $(V,Y', f_{n+1}; \mathcal{J})$ be a tuple which is 
equivalent to $(V,Y,f_{n+1}; \mathcal{J})$. Then for some $P \in {\widehat{G}}_{f_{n+1}}$, some 
$w \in V$, and every $Z \in {\widehat{\mathfrak{g}}}_{f_{n+1}}$ we have 
\begin{align}
{\ell }_{Y'}(Z) & = {\ell }_{P(Y+L_{w, f_{n+1}})P^{-1}}(Z) \notag \\
& = {\ell }_{PYP^{-1}}(Z) + {\ell }_{PL_{w, f_{n+1}}P^{-1}}(Z) \notag \\
& = {\ell}_Y(P^{-1}ZP) + {\ell }_{L_{Pw, Pf_{n+1}}}(Z), \notag \\
& \hspace{1in}\parbox[t]{2.5in}{since $\mathrm{tr}(PYP^{-1}Z) = \mathrm{tr}(Y(P^{-1}ZP)) $ \\
and ${\widehat{G}}_{f_{n+1}}\subseteq \Spp (V, \mathcal{J})$.} \notag \\
& = ({\Ad }^T_{P^{-1}}{\ell }_Y)(Z) + {\ell }_{L_{Pw, f_{n+1}}}(Z), \quad \mbox{since 
$P \in {\widehat{G}}_{f_{n+1}}$} \notag \\
& = {\Ad }^T_{P^{-1}}{\ell }_Y(Z) , \quad \mbox{since $L_{Pw, f_{n+1}} \in \widehat{\mathfrak{g}} = 
({\widehat{\mathfrak{g}}}_{f_{n+1}})^{\circ}$ and $Z\in {\widehat{\mathfrak{g}}}_{f_{n+1}}$.} \notag
\end{align}
Thus the map of cotypes to ${\widehat{G}}_{f_{n+1}}$ coadjoint orbits, induced by the 
correspondence (\ref{eq-three}), is well defined. The induced map is injective, for if 
$({\ell }_{Y'})_{|{\widehat{\mathfrak{g}}}_{f_{n+1}}}$ lies in the ${\widehat{G}}_{f_{n+1}}$ coadjoint 
orbit through $({\ell }_Y)_{|{\widehat{\mathfrak{g}}}_{f_{n+1}}}$, then for some $P \in {\widehat{G}}_{f_{n+1}}$ 
we have ${\ell }_{Y'} = {\Ad }^T_{P^{-1}}{\ell }_Y = {\ell }_{PYP^{-1}}$ on $\widehat{\mathfrak{g}}$. So 
$Y' - PYP^{-1} \in ({\widehat{\mathfrak{g}}}_{f_{n+1}})^{\circ}$. Hence for some $w \in V$ we have 
$Y' - PYP^{-1} = L_{w, f_{n+1}}$, that is, $Y' = P(Y + L_{P^{-1}w, f_{n+1}})P^{-1}$. Thus the tuples 
$(V,Y,f_{n+1}; \mathcal{J})$ and $(V, Y', f_{n+1}; \mathcal{J})$ are equivalent and thus correspond 
to the same cotype. Since every element of $({\widehat{\mathfrak{g}}}_{f_{n+1}})^{\ast}$ may be 
written as $({\ell }_Y)_{|{\widehat{\mathfrak{g}}}_{f_{n+1}}}$ for some $Y \in \widehat{\mathfrak{g}}$, the 
induced map is surjective. \hfill $\square $ \medskip

\noindent \textbf{Fact 5} {\tiny $\begin{pmatrix} 0 & 0 & 0 \\ 0 & 0 & 0 \\ \xi & 0 & 0 \end{pmatrix}$}$\in 
{\widehat{\mathfrak{g}}}_{f_{n+1}}$ is invariant under conjugation by elements of ${\widehat{G}}_{f_{n+1}}$. \medskip 

\noindent \textbf{Proof} The proof is a straightforward calculation. The details are omitted. \hfill $\square $ \medskip 

\noindent \textbf{Claim 6} The tuple $({\mathbb{R} }^{2n}, Y, f_{n+1}; \mathcal{J})$, where 
$Y=${\tiny $\begin{pmatrix} \zeta & 0 & 0 \\ d & 0 & 0 \\ \xi & (Jd)^T & -\zeta \end{pmatrix}$}$\in \widehat{\mathfrak{g}}$, 
is equivalent to the tuple $({\mathbb{R}}^{2n},\mbox{\tiny $\begin{pmatrix} 0 & 0 & 0 \\ 0 & 0 & 0 \\ \xi & 0 & 0 \end{pmatrix}$}, f_{n+1}; \mathcal{J})$. \medskip 

\noindent \textbf{Proof} $Y = L_{w',f_{n+1}}$, where $w' = \zeta e_0 + d + \onehalf \xi f_{n+1} \in 
{\mathbb{R} }^{2n+2}$. Choose $w = -\zeta e_0 -d \in {\mathbb{R} }^{2n+2}$. Then $Y +L_{w, f_{n+1}}  =${\tiny $\begin{pmatrix}
0 & 0 & 0 \\ 0 & 0 & 0 \\ \xi & 0 & 0 \end{pmatrix}$}. \hfill $\square $ \medskip 

Consider the tuple $(V,Y, f_{n+1}; \mathcal{J})$. If $V = V_1 \oplus V_2$, where $V_1$ and 
$V_2$ are $Y$-invariant, $\mathcal{J}$ perpendicular $\mathcal{J}$ nondegenerate subspaces 
of $(V, \mathcal{J})$ and $f_{n+1}\in V_1$, the cotype $\nabla $, represented by the tuple 
$(V,Y, f_{n+1}; \mathcal{J})$ is \emph{decomposable} into a cotype $\nabla $, represented by the 
tuple $(V_1, Y_{|V_1}, f_{n+1}; \mathcal{J}_{|V_1})$, and a type $\Delta $, represented by the pair 
$(V_2, Y_{|V_2}; {\mathcal{J}}_{|V_2})$. We say that $\nabla $ is the \emph{sum} of the cotype 
$\nabla $ and the type $\Delta $. If no such decomposition exists, then $\nabla $ is said to be 
\emph{indecomposable}. The pair $({\mathbb{R} }^{2n+2}, \mbox{\tiny $\begin{pmatrix} 0 & 0 & 0 \\
0 & 0 & 0 \\ \xi & 0 & 0 \end{pmatrix}$}; \mathcal{J})$ with $\xi =0$ is the indecomposable 
zero type ${\mathbf{0}}_{2n+2}$. If $\xi \ne 0$ then the cotype $\nabla $, represented by the nipotent 
tuple $({\mathbb{R}}^{2n+2}, \mbox{\tiny $\begin{pmatrix} 0 & 0 & 0 \\
0 & 0 & 0 \\ \xi & 0 & 0 \end{pmatrix}$}, f_{n+1}; \mathcal{J})$ of height $1$, is decomposable into the 
sum of the zero type ${\mathbf{0}}_{2n}$ and the indecomposable cotype ${\nabla}_1(0), \, \xi \ne 0$ 
of height $1$ and \emph{modulus} $\xi$, represented by the nilpotent tuple 
$({\mathbb{R} }^2, \mbox{\tiny $\begin{pmatrix} 0 & 0 \\ \xi & 0 \end{pmatrix}$}, f_2; {\mathcal{J}}_2 =${\tiny 
$\begin{pmatrix} 0 & -1 \\ 1 & 0 \end{pmatrix}$}$)$ of height $1$. \medskip

We explain the geometric meaning of the modulus $\xi $ in the cotype ${\nabla }_1(0),$ \linebreak  
$\xi \ne 0$. Let 
$Y_1 =${\tiny $\begin{pmatrix} 0 & 0 & 0 \\ 0 & 0 & 0 \\ 1 & 0 & 0 \end{pmatrix}$} and 
$Y_{\xi }=${\tiny $\begin{pmatrix} 0 & 0 & 0 \\ 0 & 0 & 0 \\ \xi & 0 & 0 \end{pmatrix}$}$=\xi Y_1$. We have 
\begin{align*}
{\nabla }_1(0), \, 1 & = \left\{ (V,Y, f_{n+1}; \mathcal{J} ) \setrule \, \parbox[t]{1.75in}{$Y= P(Y_1+L_{w,f_{n+1}} )P^{-1}$ \\ 
for $P \in {\widehat{G}}_{f_{n+1}}$ and $w \in V$}  \right\} .
\end{align*}
Now 
\begin{align*}
\{ (V,\widehat{Y}, f_{n+1}; \mathcal{J}) \setrule \, \widehat{Y} &= \xi Y = \xi P (Y_1 + L_{w,f_{n+1}})P^{-1} \} \\
&\hspace{-1in}= \{ (V, \widehat{Y}, f_{n+1}; \mathcal{J}) \setrule \, \widehat{Y} = 
P(Y_{\xi } + L_{\widehat{w}, f_{n+1}})P^{-1} \} = {\nabla }_1(0), \, \xi \ne 0, 
\end{align*}
where $\widehat{w} = \xi w$. Here $P \in {\widehat{G}}_{f_{n+1}}$ and $w \in V$ are those given in 
${\nabla }_1(0),\, 1$ for the tuple $(V,Y,f_{n+1}; \mathcal{J})$. The ${\widehat{G}}_{f_{n+1}}$ coadjoint 
orbit ${\mathcal{O}}_1$ through $({\ell }_{Y_1})_{|{\widehat{\mathfrak{g}}}_{f_{n+1}}}$ and the 
${\widehat{G}}_{f_{n+1}}$ coadjoint orbit ${\mathcal{O}}_{\xi}$ through 
$({\ell }_{Y_{\xi }})_{|{\widehat{\mathfrak{g}}}_{f_{n+1}}}$ are \emph{diffeomorphic} via the mapping 
${\Ad }^T_{P^{-1}}({\ell }_{Y_1})_{|{\widehat{\mathfrak{g}}}_{f_{n+1}}} \mapsto 
{\Ad }^T_{P^{-1}}({\ell }_{Y_{\xi}})_{|{\widehat{\mathfrak{g}}}_{f_{n+1}}}$, since the set of $P \in 
{\widehat{G}}_{f_{n+1}}$ defining the cotypes ${\nabla}_1(0), \, 1$ and ${\nabla }_1(0), \, \xi \ne 0$ are the 
same. The coadjoint orbits ${\mathcal{O}}_1$ and ${\mathcal{O}}_{\xi}$ are symplectic manifolds with 
their natural symplectic forms ${\omega }_{{\mathcal{O}}_1}$ and ${\omega }_{{\mathcal{O}}_{\xi }}$, respectively, 
see \cite[p.288]{cushman-bates}. Since $Y_{\xi} = \xi Y_1$ it follows that ${\ell }_{Y_{\xi}} = \xi {\ell }_{Y_1}$. 
Hence ${\omega}_{{\mathcal{O}}_{\xi}} = \xi {\omega }_{{\mathcal{O}}_1}$. Thus the symplectic manifolds 
$({\mathcal{O}}_{\xi }, {\omega }_{{\mathcal{O}}_{\xi}})$ and $({\mathcal{O}}_1, {\omega }_{{\mathcal{O}}_1})$ 
are \emph{not} symplectically diffeomorphic when $\xi \in \mathbb{R} \setminus \{ 0 , 1 \}$. The geometric 
meaning of the modulus for the cotype ${\nabla}_1(0), \, \xi \ne 0$ is: $\xi $ parametrizes a family 
$({\mathcal{O}}_{\xi}, {\omega }_{{\mathcal{O}}_{\xi}})$ of symplectic manifolds where the base manifolds 
${\mathcal{O}}_{\xi}$ are diffeomorphic, but are not pairwise symplectically diffeomorphic.  \bigskip 

\noindent {\Large \textbf{2A. Appendix. Intrinsic description}}\bigskip

Following Wallach \cite{wallach} we give an intrinsic description of the coadjoint orbits of 
the Heisenberg group $H$ on ${\mathfrak{h}}^{\ast }$, the dual of its Lie algebra $\mathfrak{h}$. \medskip

Recall that the Heisenberg group $H$ is a subset of $({\mathbb{R} }^{2n}, \omega ) \times \mathbb{R} $ with multiplication 
$(x,t) \cdot (y,s) = \big( x+y, t+s + \onehalf \omega (x,y) \big)$. $H$ is a Lie group with Lie algebra 
$\mathfrak{h}$ having the bracket $[(\xi ,t), (\eta ,s)] = (0, \omega (\xi ,\eta )]$. Its exponential map 
$\exp :\mathfrak{h} \rightarrow H:(\xi ,s) \mapsto (\xi ,s)$ is the identity map. \medskip

The first step toward computing the coadjoint action of $H$ is to compute the adjoint map. 
\begin{displaymath}
{\Ad}_{(x,t)}: \mathfrak{h} \rightarrow \mathfrak{h}: (\xi ,s) \mapsto (x,t)\cdot (\xi ,s) \cdot (x,t)^{-1} 
\end{displaymath}
for every $(x,t) \in H$. We get 
\begin{align} 
{\Ad }_{(x,t)}(\xi ,s) & = (x,t)\cdot (\xi ,s) \cdot (-x,-t) \notag \\
& = (x,t) \cdot \big( \xi -x, x-t +\onehalf \omega (\xi , -x) \big) \notag \\
& = (x+\xi -x, t+\big( x-t -\onehalf \omega (\xi ,x) + \onehalf (x, \xi -x) \big) \notag \\
& = \big( \xi , s + \omega (x, \xi ) \big) . \tag*{$\square $}
\end{align} 
Let $f \in {\mathfrak{h}}^{\ast }$. For every $(\xi ,s)\in \mathfrak{h}$ we have $f(\xi ,s) = {\lambda }_f(\xi ) + sf(0,1)$, where ${\lambda }_f \in ({\mathbb{R} }^{2n})^{\ast }$. So the coadjoint action \raisebox{1pt}{\tiny $\bullet$} of $H$ on 
${\mathfrak{h}}^{\ast}$ is 
\begin{align}
\big( (x,t) \raisebox{1pt}{\tiny $\bullet$} f \big)(\xi ,s ) & = ({\Ad}^T_{(x,t)^{-1}}f)(\xi ,s) \notag \\
& = {\lambda }_f(\xi) + \big( s - \omega (x, \xi ) \big) f(0,1) \notag \\
& = \big( {\lambda }_f(\xi ) + sf(0,1)\big) - \omega (x,\xi) f(0,1) \notag \\
& = f(\xi ,s) - \big({\omega }^{\sharp}(x)\xi \big) f(0,1). \tag*{$\square $}
\end{align}
Since ${\lambda}_f \in ({\mathbb{R} }^{2n})^{\ast}$ and $\omega $ is nondegenerate, there is $y_f \in {\mathbb{R} }^{2n}$ 
such that $y_f = {\omega }^{\flat}({\lambda }_f)$. For $\mu = f(0,1) \ne 0$ we get 
\begin{align}
\big( ({\mu }^{-1}y_f, 0) \raisebox{1pt}{\tiny $\bullet$} f \big) (x,t) & = 
f(\xi ,s) - {\omega }^{\sharp}({\mu}^{-1}y_f)(\xi ) \mu \notag \\
& = f(\xi ,s) - ({\omega }^{\sharp}({\omega }^{\flat}{\lambda }_f)) (\xi) \notag \\
& = {\lambda }_f(\xi ) + \mu s -{\lambda }_f(\xi ) = \mu s. \tag*{$\square $} 
\end{align}
For $\mu \ne 0$ let $h_{\mu }: \mathfrak{h} \rightarrow \mathbb{R} : (\xi ,s) \mapsto \mu s$. Then $h_{\mu } \in 
{\mathfrak{h}}^{\ast }$. The above calculation shows that for every $\mu \ne 0$ we have 
$({\mu }^{-1}y_f,0) \raisebox{1pt}{\tiny $\bullet$} f = h_{\mu }$. Thus we have proved \medskip 

\vspace{-.15in}\noindent \textbf{Proposition 2A.1} When $\mu \ne 0$ every $f = {\lambda }_f + h_{\mu } \in 
{\mathfrak{h}}^{\ast }$ lies in the $H$ coadjoint orbit of $h_{\mu }$. \medskip 

Now we show that \medskip 

\noindent \textbf{Proposition 2A.2} The $H$ coadjoint orbit ${\mathcal{O}}_{\mu} = 
H \raisebox{1pt}{\tiny $\bullet$} h_{\mu}$ through $h_{\mu } \in {\mathfrak{h}}^{\ast }$ is the 
symplectic manifold $({\mathbb{R} }^{2n}, \mu {\omega }_{|{\mathcal{O}}_{\mu}})$. \medskip 

\noindent \textbf{Proof} First we compute the isotropy group $H_{h_{\mu}} = \{ h \in H \setrule \, 
h \raisebox{1pt}{\tiny $\bullet$} h_{\mu} = h_{\mu} \} $ of the $H$ coadjoint action at $h_{\mu }$. Since 
$((x,t) \raisebox{1pt}{\tiny $\bullet$} h_{\mu })(\xi ,s) = h_{\mu }(\xi ,s) - \mu {\omega }^{\sharp}(x)\xi$, 
we see that $(\xi ,s) \in H_{h_{\mu }}$ implies that $0 = {\mu }^{-1}{\omega }^{\sharp}(x)\xi$ 
for every $(\xi ,s) \in \mathfrak{h}$. Hence $x=0$, since $\mu \ne0$ and $\omega $ is nondegenerate. 
Thus $H_{h_{\mu }}$ is the center $Z = \{ (0,t) \in H \setrule \, t \in \mathbb{R} \}$ of $H$. By definition of the coadjoint 
orbit the map $\psi : H/Z \rightarrow {\mathcal{O}}_{\mu }: hZ \mapsto h \raisebox{1pt}{\tiny $\bullet$} h_{\mu }$ 
is smooth and bijective, as is the map $\eta : {\mathbb{R} }^{2n} \rightarrow {\mathcal{O}}_{\mu }: x \mapsto 
(x,0) \raisebox{1pt}{\tiny $\bullet$} h_{\mu}$. Since the map 
$\theta : H/Z \rightarrow {\mathbb{R} }^{2n} \times \{ 0 \} : (x,t)Z \mapsto (x,0)$ is smooth and bijective, it follows 
that the coadjoint orbit ${\mathcal{O}}_{\mu }$ is diffeomorphic to ${\mathbb{R} }^{2n}$. If $\alpha = 
(\xi ,0)$ and $\beta = (\eta ,0)$ lie in $\mathfrak{h}$, for $\pi \in {\mathcal{O}}_{\mu }$ 
\begin{displaymath}
{\omega }_{{\mathcal{O}}_{\mu }}(\pi )(X^a(\pi), X^b(\pi) ) = \pi ([a,b]) = \pi (0, \omega (\xi ,\eta )) = 
\mu {\omega }_{|{\mathcal{O}}_{\mu }}(\pi). 
\end{displaymath}
Here $X^c(\pi) = (\smalldbydt \hspace{-10pt}\exp tc) \raisebox{1pt}{\tiny $\bullet$} \pi $ for $c \in \mathfrak{h}$.

\section{The momentum mapping of ${\mathrm{H}}_{2n+1}$}

In this section we show that the ${\mathbb{R} }^{2n}$ cocycle of the momentum map $\Phi $ (\ref{eq-one}) 
of the affine action of ${\mathbb{R} }^{2n}$ on $({\mathbb{R} }^{2n}, \omega )$ becomes the value of a modulus for 
a coadjoint orbit of the Heisenberg group ${\mathrm{H}}_{2n+1}$ coming from the action ${\Phi }^{\vvee}$ 
(\ref{eq-two}) on $({\mathbb{R} }^{2n}, \omega )$. \medskip 

We determine that image of the momentum map $J: {\mathbb{R} }^{2n} \rightarrow {\mathfrak{h}}^{\ast }$ of 
the Hamiltonian action ${\Phi }^{\vvee}$. By definition $J(v)X^{\vvee} = J^{X^{\vvee}}(v)$, where 
$J^{X^{\vvee}}: {\mathbb{R} }^{2n} \rightarrow \mathbb{R} : v \mapsto \omega (x,v) + \xi $ for $X^{\vvee}=${\tiny 
$\begin{pmatrix} 0 &0 &0 \\ \widetilde{x} & 0 & 0 \\ \xi & (J\widetilde{x})^T & 0 \end{pmatrix}$}. The action 
${\Phi }^{\vvee}$ is transitive. Thus 
\begin{displaymath}
J({\mathbb{R} }^{2n}) = J({\Phi }^{\vvee}_{{\mathbb{R}}^{2n}}(0)) = {\Ad }^T_{{\mathbb{R} }^{2n}}(J(0)) = 
H \raisebox{1pt}{\tiny $\bullet$} J(0), 
\end{displaymath}
the $H$ coadjoint orbit through $J(0)$. But $J(0)X^{\vvee} = J^{X^{\vvee}}(0) = \xi $. So 
$J(0) = E^{\ast }_{0, 2n+1}$, since 
\begin{displaymath}
J(0)\mbox{\footnotesize $\begin{pmatrix} 0 & 0 & 0 \\ \widetilde{x} & 0 & 0 \\ 
\xi & {\omega }^{\sharp }(\widetilde{x} )& 0 \end{pmatrix} $} = \xi = E^{\ast }_{0, 2n+1}\mbox{\footnotesize $\begin{pmatrix} 0 & 0 & 0 \\ \widetilde{x} & 0 & 0 \\ \xi & {\omega }^{\sharp}(\widetilde{x}) & 0 \end{pmatrix} $} .
\end{displaymath}
Thus the $H$ coadjoint orbit ${\mathcal{O}}_1$ through ${\ell }_{E_{0, 2n+1}} = E^{\ast}_{0, 2n+1} \in 
{\mathfrak{h}}^{\ast}$ corresponds to the cotype represented by the tuple $({\mathbb{R} }^{2n+2}, E_{0, 2n+1}, 
f_{n+1}; \mathcal{J})$. This cotype is the sum of the indecomposable cotype ${\nabla }_1(0), \, 1$ with modulus $1$ 
and the zero type ${\mathbf{0}}_n$.

\section{Representation of $\mathfrak{h}$ corresponding to ${\mathcal{O}}_1$}

According to the theory of Kirillov \cite{kirillov}, associated to the coadjoint orbit ${\mathcal{O}}_1$ of the 
Heisenberg group $H$ through ${\ell }_X \in {\mathfrak{h}}^{\ast }$, where $X=${\tiny 
$\begin{pmatrix} 0 & 0 & 0 \\ 0 & 0 & 0 \\ 1 & 0 & 0 \end{pmatrix}$}$\in \mathfrak{h}$, there is 
an irreducible unitary representation of $\mathfrak{h}$ by skew Hermitian differential operators on 
$C^{\infty}({\mathbb{R} }^n, \C)$. We find this representation using geometric quantization, see Kostant \cite{kostant} 
or \'{S}niatycki \cite{sniatycki}. \medskip

Let ${\mathbb{R} }^{2n} = T^{\ast }{\mathbb{R} }^n$ be the cotangent bundle of ${\mathbb{R} }^n$ with coordinates $(x,y)$ and 
standard symplectic form $\omega \big( (x,y), (z,w) \big) = \langle x, w \rangle - \langle y, z \rangle $. 
Here $\langle \, \, , \, \, \rangle $ is the Euclidean inner product on ${\mathbb{R} }^n$. Consider the trivial 
complex line bundle $\rho : L = {\mathbb{R} }^{2n} \times \C \rightarrow {\mathbb{R} }^{2n}:(x,y,z) \mapsto (x,y)$. 
Let $a$, $b \in {\rho }^{-1}(x,y)$. 
Then $h(x,y)(a,b) = a\overline{b}$ defines a Hermitian inner product on $L$. A smooth section 
$\widetilde{\sigma }: {\mathbb{R} }^{2n} \rightarrow L: (x,y) \mapsto \big( (x,y), \sigma (x,y) \big)$ of $L$ will be 
identified with the smooth complex valued function $\sigma :{\mathbb{R} }^{2n} \rightarrow \C$. Let $s \in 
C^{\infty}({\mathbb{R} }^{2n}, \C)$ such that $s(x,y)=1$ for every $(x,y) \in {\mathbb{R} }^{2n}$. Then $s$ is a smooth 
section of $L$. We have a $1$-form $\theta = \langle y, \dee x \rangle $ on $T^{\ast }{\mathbb{R} }^n$, which gives 
rise to the symplectic form $\omega = - \dee \theta $. For every smooth vector field $X$ on 
$T^{\ast }{\mathbb{R} }^n$ and every smooth section $\widetilde{\sigma }$ of $L$ let 
\begin{displaymath}
{\nabla }_X\widetilde{\sigma }(x,y) = 2\pi \mathrm{i}(X \lefthook \theta )(x,y)\widetilde{\sigma}(x,y), \quad 
\mbox{for every $(x,y)\in T^{\ast }{\mathbb{R} }^n$.}
\end{displaymath}
For every $f \in C^{\infty}({\mathbb{R} }^{2n}, \C) $ the expression 
${\nabla }_X(f\widetilde{\sigma }) = (L_Xf)\cdot \widetilde{\sigma} + f \cdot ({\nabla }_X\widetilde{\sigma})$ defines a \emph{connection} on the space of smooth sections $\Gamma (L)$ of the line bundle $L$. Using the section $s(x,y)=1$ and the definition of ${\nabla }_Xs$ we obtain 
\begin{displaymath}
{\nabla}_Xf = L_Xf + 2\pi \mathrm{i} (X \lefthook \theta )f \quad \mbox{for every $f \in C^{\infty}(L, \C)$.} 
\end{displaymath}
Then for $j=1, \ldots , n$ 
\begin{displaymath}
{\nabla}_{\frac{\partial }{\partial x_j}}f = \frac{\partial f}{\partial x_j} + 2\pi \mathrm{i}y_j \, f \quad 
\mathrm{and} \quad {\nabla}_{\frac{\partial }{\partial y_j}}f = \frac{\partial f}{\partial y_j} . 
\end{displaymath}
If $F\in C^{\infty}({\mathbb{R} }^{2n}, \mathbb{R}) $ the corresponding Hamiltonian vector field on $(T^{\ast }{\mathbb{R} }^n, \omega )$ is 
$X_F = \langle \frac{\partial F}{\partial y}, \frac{\partial }{\partial x} \rangle - 
\langle \frac{\partial F}{\partial x}, \frac{\partial }{\partial y} \rangle $. \medskip 

From geometric quantization we obtain the prequantization operator $\mathcal{P}$ on $\Gamma (L)$ given by 
\begin{align*}
\mathcal{P}(F) & = -{\nabla }_{X_F} + 2\pi \mathrm{i}\ F \\
& = -\langle \frac{\partial F}{\partial y}, \frac{\partial}{\partial x} +2\pi \mathrm{i} \, y \rangle + 
\langle \frac{\partial F}{\partial x}, \frac{\partial }{\partial y} \rangle + 2\pi \mathrm{i}\, F \\
& = -X_F -2\pi \mathrm{i}\, \big( \langle y, \frac{\partial F}{\partial y} \rangle - F \big) .
\end{align*}

Consider the Lie algebra mapping 
\begin{displaymath}
\mathfrak{J}: (\mathfrak{h}, [\, \, , \, \, ]) \rightarrow (C^{\infty}({\mathbb{R} }^{2n}, \C ), \{ \, \, , \, \, \} ): (\xi ,\eta ,s) 
\mapsto J^{(\xi ,\eta , s)}, 
\end{displaymath}
where 
\begin{displaymath}
J^{(\xi ,\eta , s)}(x,y) = \omega \big( (\xi ,\eta ), (x,y) \big) +s = \langle \xi , y \rangle -\langle \eta ,x \rangle +s. 
\end{displaymath}
Here $J:{\mathbb{R} }^{2n} \rightarrow {\mathfrak{h}}^{\ast }$, where $J(x,y)(\xi ,\eta ,s) = J^{(\xi ,\eta ,s)}(x,y)$ is the momentum map of the Heisenberg group $H$ acting on $({\mathbb{R} }^{2n}, \omega )$. The mapping 
$\mathfrak{J}$ is a homomorphism of Lie algebras, namely, $J^{[(\xi ,\eta , s), ({\xi }' , {\eta }', s')]} = 
\{ J^{(\xi , \eta , s)}, J^{({\xi }', {\eta }', s')} \}$. Now 
\begin{align*}
\mathcal{P}(J^{(\xi ,\eta ,s)}) & = - X_{J^{(\xi ,\eta , s)}} - 2\pi \mathrm{i}\big( \langle y, 
\frac{\partial J^{(\xi ,\eta ,s)}}{\partial y} \rangle - J^{(\xi , \eta ,s)} \big) \\
& = -\langle \xi , \frac{\partial }{\partial x} \rangle + \langle \eta , \frac{\partial }{\partial y} \rangle - 
2\pi \mathrm{i}\big( \langle \eta , x \rangle - s \big)
\end{align*}
is the prequantization operator, which gives a representation of $\mathfrak{h}$ on 
$C^{\infty}({\mathbb{R} }^{2n},$ $ \C)$, namely, $(\xi , \eta , s) \mapsto \mathcal{P}(J^{(\xi ,\eta ,s)})$. This completes prequantization. \medskip

Using the polarization $\frac{\partial }{\partial y} =0$ of the coadjoint orbit ${\mathcal{O}}_1$, which is 
${\mathbb{R} }^{2n}$, we obtain the quantization operator 
\begin{equation}
\mathcal{Q}(\xi ,\eta ,s) = -\langle \xi , \frac{\partial }{\partial x} \rangle + 2\pi \mathrm{i}\, (s - \langle \eta , x 
\rangle ) 
\label{eq-six}
\end{equation}
on $C^{\infty}({\mathbb{R}}^n , \C)$, where ${\mathbb{R} }^n$ has coordinate $x$. On $C^{\infty}({\mathbb{R} }^{2n}, \C)$ place 
the Hermitian inner product $\langle f, g \rangle = \int_{{\mathbb{R} }^n} f \cdot g$. The quantization operator 
$\mathcal{Q}(\xi ,\eta , s)$ is skew Hermitian, that is, $\langle \mathcal{Q}(\xi , \eta ,s)f ,g \rangle = 
-\langle f, \mathcal{Q}(\xi , \eta , s)g \rangle $. Using the inner product $\langle \, \, , \, \, \rangle $, we  
complete $C^{\infty}_c({\mathbb{R} }^n, \C)$, the space of smooth complex valued functions on ${\mathbb{R} }^n$ with 
compact support, to the Hilbert space $(\mathcal{H}, \langle \, \, ,  \, \rangle )$. Then 
\begin{displaymath}
\mathfrak{h} \rightarrow \gl (\mathcal{H}, \langle \, \, , \, \, \rangle ): (\xi , \eta ,s) \mapsto 
\mathcal{Q}(\xi , \eta , s)
\end{displaymath}
is the infinitesimalization of an irreducible unitary representation of the Heisenberg group $H$. \medskip 

We now find the irreducible unitary representation of the Heisenberg group $H$, whose infinitesimalization 
is given by the quantum operator $\mathcal{Q}$ (\ref{eq-six}). Consider the representation ${\widetilde{S}}_1 = S_1 \comp \psi $, where 
\begin{displaymath}
S_1: H \rightarrow \mathrm{U}(\mathcal{H}, \langle \, \, , \, \, \rangle ): (x',y',t') \mapsto 
\big( f \mapsto S_1(x',y',t')f \big) 
\end{displaymath}
and $\big( S_1(x',y',t')f \big) (z) = {\mathrm{e}}^{2\pi \mathrm{i} [t' + \langle x', z + \onehalf y' \rangle ]}f(z-y')$ of $H$, see Wallach \cite[p.107]{wallach}. Here 
$\psi $ is the group homomorphism 
\begin{displaymath}
\psi : H \rightarrow H: (x,y,t) \mapsto (x',y',t') = (-y,x,t). 
\end{displaymath}
Let $(\xi ,\eta ,s) \in \mathfrak{h}$. Then 
\begin{displaymath}
\exp (-u\eta , u \xi ,us) = (-u\eta , u \xi , us) = (x',y',t') \in H.
\end{displaymath}
So the infinitesimalization of ${\widetilde{S}}_1$ is %
\begin{align}
\mbox{${\displaystyle \frac{d}{du}}
\rule[-10pt]{.5pt}{25pt} \raisebox{-10pt}{$\, {\scriptstyle u=0}$}$} \big( {\widetilde{S}}_1(u\xi, u\eta ,us)f \big) (z) & = 
\mbox{${\displaystyle \frac{d}{du}}
\rule[-10pt]{.5pt}{25pt} \raisebox{-10pt}{$\, {\scriptstyle u=0}$}$} 
{\mathrm{e}}^{2\pi \mathrm{i} [us - \langle u\eta , z + \onehalf u\xi \rangle ]}f(z-u\xi ) \notag \\
& = 2\pi \mathrm{i}[s - \langle \eta ,z \rangle ]f(z) - \langle \xi , \frac{\partial f}{\partial z} \rangle \notag \\
& = (\mathcal{Q}(\xi , \eta , s)f)(z),  \tag*{$\square $}
\end{align}
as desired. \medskip 

An irreducible unitary representation of $H$ corresponding to the coadjoint orbit 
${\mathcal{O}}_{\xi}$ is $S_{\xi}$, where 
\begin{displaymath}
(S_{\xi }(x,y,t)f)(z) = {\mathrm{e}}^{2\pi \mathrm{i} \xi [t + \langle -y, z + \onehalf x \rangle ]}f(z-x) . 
\end{displaymath}
Here the modulus $\xi $ of the coadjoint orbit ${\mathcal{O}}_{\xi }$ is a parameter in the irreducible 
unitary representation. Because every coadjoint orbit of the Heisenberg group $H$ is of the form 
${\mathcal{O}}_{\xi}$ for some value of the real nonzero parameter $\xi $, we have found 
all the irreducible unitary representations of the Heisenberg group up to equivalence.


\begin{thebibliography}{9}

\bibitem{cushman-bates} R.H. Cushman and L.M. Bates, \lq \lq Global aspects of classical 
integrable systems\rq \rq , second edition, Birkh\"{a}user, Basel, 2015.

\bibitem{cushman-vanderKallen} R.Cushman and W. van der Kallen, Adjoint and coadjoint orbits of the 
Poincar\' {e} group, \textit{Acta Appl. Math.} \textbf{90} (2006) 65--89. revised version: arxiv:math. math/0305442v3.

\bibitem{kirillov} A.A. Kirillov, \lq \lq Lectures on the orbit method\rq \rq , Graduate studies in 
mathematics \textbf{64} American Mathematical Society, Providence, R.I., 2004. 

\bibitem{kostant} B. Kostant, Quantization and unitary representations, part I, in: 
\lq \lq Lectures in modern analysis and applications III, ed. C.T. Taam, 
\textit{Lecture Notes in Mathematics} \textbf{170} 87--208, Springer Verlag, New York, 1970. 

\bibitem{sniatycki} J. \'{S}niatycki, \lq \lq Geometric quanization and quantum mechanics\rq \rq , 
Applied mathematics series \textbf{36} Springer Verlag, New York, 1980.

\bibitem{wallach} N. Wallach, \lq \lq Symplectic geometry and Fourier analysis\rq \rq, second edition, 
Dover, Mineola, N.Y., 2018.

\end{thebibliography}
\end{document}